\documentclass[12pt]{article} 
\newcommand{\be}{\begin{equation}}
\newcommand{\ee}{\end{equation}}
\newcommand{\ba}{\begin{array}}
\newcommand{\ea}{\end{array}}

\newcommand{\dvg}{{\rm div\;}}
\newcommand{\grad}{{\rm \; grad\;}}
\newcommand{\curl}{{\rm \; curl\;}}

\newcommand{\RR}{\mathbb R}

\newcommand{\bea}{\begin{eqnarray*}}
\newcommand{\eea}{\end{eqnarray*}}
\newcommand{\bean}{\begin{eqnarray}}
\newcommand{\eean}{\end{eqnarray}}
\newcommand{\proof}{\vspace{1ex}\noindent{\em Proof}. \ }
\def\ds{\displaystyle}
\def\nm{\noalign{\medskip}}

\newcommand\eps{\varepsilon}

\newtheorem{lemma}{Lemma}[section]

\newtheorem{theorem}{Theorem}[section]

\newtheorem{proposition}{Proposition}[section]
\usepackage{amsmath, amssymb}
\usepackage{latexsym}
\newcommand{\R}{\mathbb{R}}

\def\Box{\leavevmode\vbox{\hrule
     \hbox{\vrule\kern5pt\vbox{\kern5pt}%
           \vrule}\hrule}}
\renewcommand{\square}{\hfill$\Box$}
\begin{document}

\title{Reconstruction of closely spaced small
 inhomogeneities via boundary measurements for the full time-dependent Maxwell's
  equations}
\author{Christian Daveau \thanks{
D\'epartement de Math\'ematiques, Site Saint-Martin II, BP 222, \&
Universit\'e de Cergy-Pontoise, 95302 Cergy-Pontoise Cedex, France
(Email: christian.daveau@math.u-cergy.fr).}   Abdessatar Khelifi
\thanks{ D\'epartement de Math\'ematiques \& Informatique
Facult\'e des Sciences, 7021 Zarzouna - Bizerte, Tunisia
(Email:abdessatar.khelifi@fsb.rnu.tn).} Anton Sushchenko\thanks{
ETIS \& UMR CNRS 8051, 6 avenue du Ponceau, BP 44, 95014
Cergy-Pontoise Cedex, France (Email: anton.sushchenko@ensea.fr).}}
\maketitle \abstract{ We consider for the full time-dependent
Maxwell's equations the inverse problem of identifying locations
 and certain properties of small electromagnetic inhomogeneities
in a homogeneous background medium from dynamic boundary
measurements on the boundary for a finite time interval.}\\



\section{Introduction}
The ultimate objective of the work described in this paper is to
determine locations and certain properties of the shapes of small
electromagnetic inhomogeneities in a homogeneous background medium
from dynamic boundary measurements on part of the boundary and for
finite interval in time. Using as weights particular background
solutions constructed by a geometrical control method we develop
an asymptotic method based on appropriate averaging of the partial
dynamic boundary measurements.\\
For stationary Maxwell's equations it has been known that the
Dirichlet to Neumann map uniquely determines (smooth) isotropic
electromagnetic parameters, see \cite{M2}, \cite{RK}, \cite{SIC}.
We will provide in this paper a rigorous derivation of the inverse
Fourier transform of a linear combination of derivatives of point
masses, located at the positions $z_j$ of the inhomogeneities, as
the leading order term of an appropriate averaging of (partial)
dynamic boundary measurements of the tangential components of
electric fields on part of the boundary. Our formulas may be used
to determine properties (location, relative size ) of the small
inhomogeneities in case a single, or a few (tangential) boundary
electric fields are known. Our approach differs from \cite{A},
\cite{A-al}, \cite{AMV}, \cite{AVV}, \cite{Volkov} and is expected
to lead to
very effective computational identification algorithms.\\
Our main result is given by: \\ 
{\bf Theorem 4.1\/} {\em
Let $\eta \in \RR^d$. Let $E_\alpha$ be the unique solution in $
{\cal C}^0(0, T; X(\Omega)) \cap {\cal C}^1(0, T; L^2(\Omega))$ to
the Maxwell's equations (\ref{alphwam}) with $ \varphi(x) =
\eta^{\perp} e^{i \eta \cdot x},$ $ \psi(x) = - i \sqrt{\mu_0}|
\eta | \eta^{\perp} e^{i \eta \cdot x},$ and $ f(x, t) =
\eta^{\perp} e^{i \eta \cdot x  -  i \sqrt{\mu_0}| \eta | t}.$
Suppose that $\Gamma$ and $T$ geometrically control $\Omega$, then
we have $$
\begin{array}{l}
\ds \int_0^T \int_\Gamma \Big[ \theta_\eta \cdot (\curl E_\alpha
\times {\bf n} - \curl E \times {\bf n}) +
\partial_t \theta_\eta \cdot \partial_t (\curl E_\alpha \times {\bf n}
- \curl E \times {\bf n})  \Big]~d\sigma(x)dt =  \\ \nm \ds
\alpha^2 \sum_{j=1}^m (\mu_0 - \mu_j) e^{2 i \eta \cdot z_j}
M_j(\eta) \cdot \eta \;
  + O(\alpha^2),
  \end{array}
$$ where $\theta_\eta$ is the unique solution to the Volterra
equation (\ref{eq4p}) with $g_\eta$ defined as the boundary
control  in (\ref{wetam}) and $M_j$ is the polarization tensor of
$B_j$, defined by $$ \ds (M_j)_{k, l} = e_k \cdot
(\int_{\partial B_j}
 (\nu_j +
(\frac{\mu_j}{\mu_0} - 1) \frac{\partial \Phi_j}{\partial
\nu_j}|_+ (y)) y \cdot e_l \; ds_j(y)). $$ Here $(e_1, e_2)$ is
an orthonormal basis of $\RR^d$. The term $O(\alpha^2)$ is
independent of the points $\{z_j,\quad j=1,\cdots, m\}$.
\/} \\
For discussions on closely related (stationary) identification
problems we refer the reader to \cite{SU},\cite{VV}, \cite{CMV},
and \cite{FV}.
\section{Problem formulation} Let $\Omega$ be a bounded $C^2$-domain
in $\mathbb{R}^d$, $d= 2,3$. Assume  that $\Omega$ contains a
finite number of inhomogeneities, each of the form $z_j + \alpha
B_j$, where $B_j \subset \R^d$ is a bounded, smooth domain
containing the origin. The total collection of inhomogeneities is
 $ \ds {\cal B}_\alpha
 = \ds \cup_{j=1}^{m} (z_j  + \alpha B_j)$.
  The points $z_j \in \Omega, j=1, \ldots, m,$ which determine the
  location of the inhomogeneities, are assumed to satisfy the
  following inequalities:
\be \label{f1}
 | z_j  - z_l | \geq c_0 > 0, \forall \; j \neq l  \quad
\mbox{ and } \mbox{ dist} (z_j, \partial \Omega) \geq c_0
> 0, \forall \; j.
\ee Assume that $\alpha
>0$, the common order of magnitude of the diameters of the
inhomogeneities, is sufficiently small, that these
 inhomogeneities are disjoint,  and that
their distance to $\R^d \setminus \overline{\Omega}$ is larger
than $c_0/2$.  Let $\mu_0$ and $\eps_0$ denote the permeability
and the permittivity of the background medium, and assume that
{\it $\mu_0>0$ and $\eps_0>0$ are positive constants}.  Let
$\mu_j>0$ and $\eps_j>0$ denote the permeability and the
permittivity of the j-th inhomogeneity, $z_j+\alpha B_j$, these
are also assumed to be positive constants. Introduce the
piecewise-constant magnetic permeability
\begin{equation}
\mu_\alpha(x)=\left \{ \begin{array}{*{2}{l}}
 \mu_0,\;\;& x \in \Omega \setminus \bar {\cal B}_\alpha,  \\
 \mu_j,\;\;& x \in z_j+\alpha B_j, \;j=1 \ldots m.
\end{array}
\right . \label{murhodef}
\end{equation}
If we allow the degenerate case $\alpha =0$, then the function
$\mu_0(x)$ equals the constant $\mu_0$. The electric permittivity
is defined by $\eps_\alpha(x)=\eps_0$, for all $x\in \Omega$. Let
${\bf n}={\bf n}(x)$ denote the outward unit normal vector to
$\Omega$ at a point on $\partial \Omega$, $\ds\partial_tu=
\frac{\partial u}{\partial t}$ and $\Delta$ means the Laplace operator defined by $\Delta u=
\ds\sum_{i=1}^{d}\frac{\partial^2 u}{\partial x_{i}^{2}}$. \\

In this paper, we will denote by bold letters the functional
spaces for the vector fields. Thus $H^{s}(\Omega)$ denotes the
usual Sobolev space on
 $\Omega$ and ${\bf H}^s(\Omega)$ denotes $(H^s(\Omega))^d$ and ${\bf L}^2(\Omega)$
  denotes $(L^2(\Omega))^d$. As usual for Maxwell
 equations, we need spaces of fields with square integrable curls:
\[
{\bf H}(\curl; \Omega) = \{ u \in {\bf L}^2(\Omega), \curl u \in
{\bf L}^2(\Omega) \},\] and with square integrable divergences
\[
{\bf H}(\dvg; \Omega) = \{ u \in {\bf L}^2(\Omega), \dvg u \in
{\bf L}^2(\Omega) \}.\] We will also need the following functional
spaces:
\[
Y(\Omega) = \{ u \in {\bf L}^2(\Omega), \dvg u = 0 {\rm \; in\;}
\Omega
 \},\quad X(\Omega) = {\bf H}^1(\Omega)\cap Y(\Omega),
\] and $TL^2(\partial \Omega)$ the space of vector fields on
$\partial \Omega$
  that lie in ${\bf L}^2(\partial\Omega)$. Finally, the "minimal" choice for
   the electric variational space would be
   \[
 X_N(\Omega)=\{v\in {\bf H}(\curl; \Omega)\cap {\bf H}(\dvg; \Omega) ;\quad
 v\times {\bf n}=0\quad \mbox{on } \partial\Omega\}. \]

Now, we introduce the following time-dependent Maxwell equations
(associated to the electric field) \be \label{alphwam} \left\{
\begin{array}{l}
\ds (\eps_\alpha\partial_t^2 + \curl \frac{1}{\mu_\alpha} \curl )
E_\alpha = 0 \quad
{\rm in}\; \Omega \times (0, T),\\
\nm \dvg (\eps_\alpha E_\alpha) = 0 \quad
{\rm in}\; \Omega \times (0, T),\\
\nm E_\alpha |_{t=0} = \varphi,
\partial_t E_\alpha |_{t=0} = \psi \quad {\rm in}\; \Omega,\\ \nm
E_\alpha \times {\bf n} |_{\partial \Omega \times (0, T)} = f,
\end{array}
\right. \ee where $E_\alpha \in \RR^d$ is the electric field, $f$
the boundary condition for $E_\alpha \times {\bf n}$, and
$\varphi$ and
$\psi$ the initial data.\\

Let $E$  be the solution of the Maxwell's equations in the
homogeneous domain: \be \label{wom} \left\{
\begin{array}{l}
\ds
 (\eps_0\partial_t^2 + \curl \frac{1}{\mu_0}\curl ) E = 0 \quad {\rm in}\;
\Omega \times (0, T),\\
\nm \dvg (\eps_0 E)  = 0 \quad
{\rm in}\; \Omega \times (0, T),\\
\nm E |_{t=0} = \varphi, \partial_t E|_{t=0} = \psi \quad {\rm
in}\;
\Omega,\\
\nm E \times {\bf n} |_{\partial \Omega \times (0, T)} = f.
\end{array}
\right. \ee Here $T
>0$ is a final observation time and
 $\varphi, \psi \in {\cal
C}^{\infty}(\overline{\Omega})$
 and $f \in {\cal C}^{\infty}(0, T; {\cal C}^{\infty}(\partial
 \Omega))$ are subject to the compatibility conditions
\[ \partial_{t}^{2l} f |_{t = 0} =  (\Delta^l \varphi) \times {\bf n}
|_{\partial \Omega} {\rm \; and \;} \partial_t^{2l + 1}  f |_{t
=0} =  (\Delta^l \psi ) \times {\bf n} |_{\partial \Omega}, \quad
l=1, 2, \ldots \]  it follows  that (\ref{wom}) has a unique
solution $ E \in {\cal C}^\infty ([0, T] \times
\overline{\Omega})$. It is also known (see for example \cite{NI})
that since $\Omega$ is smooth (${\cal C}^2-$ regularity would be
sufficient) the non homogeneous  Maxwell's equations
(\ref{alphwam}) have a unique weak solution $E_\alpha \in {\cal
C}^0(0, T; X(\Omega)) \cap {\cal C}^1(0, T; {\bf L}^2(\Omega))$.
Indeed, $\curl E_\alpha$ belongs to $ {\cal C}^0(0, T; X(\Omega))
\cap {\cal C}^1(0, T; {\bf L}^2(\Omega))$.

\section{Asymptotic formula}
We start the derivation of the asymptotic formula for $\ds \curl
E_\alpha\times {\bf n}$ with the following estimate.

\begin{lemma}\label{est-lem1}
The following estimate as $\alpha \rightarrow 0$ holds: \be
\label{estenergy0} || \partial_t (E_\alpha - E) ||_{L^{\infty}(0,
T; {\bf L}^2(\Omega))} + || E_\alpha - E ||_{L^{\infty}(0, T;
X_N(\Omega))}  \leq C \alpha, \ee
 where the constant $C$ is independent of
$\alpha$ and the set of points $\{  z_j\}_{j=1}^m$ provided that
assumption (\ref{f1}) holds.
\end{lemma}
\proof From (\ref{alphwam})-(\ref{wom}), it is obvious that
$E_\alpha - E \in X_N(\Omega)$, then due to the Green formula we
have for any ${\bf v} \in X_N(\Omega)$: \be \label{eq2} \ds
\int_{\Omega} \eps_0
\partial_t^2 (E_\alpha -E) \cdot{\bf v}~dx + \ds\int_{\Omega}
\frac{1}{\mu_\alpha} \curl (E_\alpha -E) \cdot \curl {\bf v}~dx
=\ee
\[\sum_{j=1}^m (\frac{1}{\mu_0} - \frac{1}{\mu_j}) \int_{z_j +
\alpha B_j} \curl E \cdot \curl {\bf v}~dx.\]  Let ${\bf
v}_\alpha$ be defined by \be\label{rel-lem1-1} \left\{
\begin{array}{l}
{\bf v}_\alpha \in X_N(\Omega),\\
\nm \curl \frac{1}{\mu_\alpha} \curl {\bf v}_\alpha = \partial_t
(E_\alpha - E) \quad {\rm in\;} \Omega.
\end{array}
\right. \ee Then,
\[
\ds  \int_{\Omega} \frac{1}{\mu_\alpha} \curl (E_\alpha -E) \cdot
\curl {\bf v}_\alpha ~dx = - \int_{\Omega} \partial_t (E_\alpha -
E) \cdot(E_\alpha -E)~dx = \] \[- \frac{1}{2} \partial_t
\int_{\Omega} |E_\alpha - E|^2~dx
\]
and by Green formula, relation (\ref{rel-lem1-1}) gives:
\[
\begin{array}{lll}
\ds \int_{\Omega} \partial_t^2 (E_\alpha -E)\cdot {\bf
v}_\alpha~dx &=& \ds \int_{\Omega}
\curl \frac{1}{\mu_\alpha} \curl \partial_t {\bf v}_\alpha \cdot {\bf v}_\alpha~dx\\
\nm &=& \ds - \int_{\Omega}
\frac{1}{\mu_\alpha} \curl \partial_t {\bf v}_\alpha \cdot \curl {\bf v}_\alpha~dx\\
\nm &=& \ds - \frac{1}{2} \partial_t \int_{\Omega}
\frac{1}{\mu_\alpha} |\curl {\bf v}_\alpha|^2~dx.\end{array}
\]
Thus, it follows from (\ref{eq2}) that
\[
\ds \eps_0\partial_t \int_{\Omega} \frac{1}{\mu_\alpha} |\curl
{\bf v}_\alpha|^2~dx +
\partial_t \int_{\Omega} |E_\alpha - E|^2~dx =\]
\[
 -2 \sum_{j=1}^m (\frac{1}{\mu_0} - \frac{1}{\mu_j}) \int_{z_j
 + \alpha B_j} \curl E \cdot \curl {\bf v}_\alpha~dx.
\]
Next,
\[
\ds | \sum_{j=1}^m (\frac{1}{\mu_0} - \frac{1}{\mu_j})
 \int_{z_j + \alpha B_j} \curl E \cdot \curl {\bf v}_\alpha | \leq C || \curl E
  ||_{{\bf L}^2({\cal B}_\alpha)} || \curl {\bf v}_\alpha
||_{{\bf L}^2(\Omega)}.
\]
Since $E \in {\cal C}^{\infty}([0, T] \times \overline{\Omega})$
we have
\[
|| \curl E ||_{{\bf L}^2({\cal B}_\alpha)} \leq || \curl E
||_{L^\infty({\cal B}_\alpha)} \alpha (\sum_{j=1}^m
|B_j|)^{\frac{1}{2}} \leq C \alpha,
\]
which gives
\[
\ds | \sum_{j=1}^m (\frac{1}{\mu_0}- \frac{1}{\mu_j}) \int_{z_j +
\alpha B_j} \curl E\cdot \curl {\bf v}_\alpha dx| \leq C \alpha ||
\curl {\bf v}_\alpha ||_{{\bf L}^2(\Omega)}
\]
and so, \be\label{rel-lem1-2} \ds \eps_0\partial_t \int_{\Omega}
\frac{1}{\mu_\alpha} |\curl {\bf v}_\alpha|^2~dx +
\partial_t \int_{\Omega} |E_\alpha - E|^2~dx \leq C \alpha
(\int_{\Omega} \frac{1}{\mu_\alpha} |\curl {\bf v}_\alpha|^2~dx +
\int_{\Omega} |E_\alpha - E|^2~dx )^{1/2}. \ee From the Gronwall
Lemma it follows that \be \label{laast1} (\int_{\Omega}
\frac{1}{\mu_\alpha} |\curl {\bf v}_\alpha|^2~dx)^{1/2} +
(\int_{\Omega} |E_\alpha - E|^2~dx)^{1/2} \leq C \alpha. \ee
Combining this last estimate (\ref{laast1}) with the fact that
\[\ds
|| \partial_t (E_\alpha - E)||_{L^{\infty}(0, T; H^{-1}(\Omega))}
\leq C || \curl {\bf v}_\alpha ||_{L^{\infty}(0, T; {\bf
L}^2(\Omega))}
\]
the following estimate holds \be \label{estenergy1} \ds ||
E_\alpha - E_0||_{L^{\infty}(0, T; {\bf L}^2(\Omega))} + ||
\partial_t (E_\alpha - E_0)||_{L^{\infty}(0, T; {\bf L}^2(\Omega))}
\leq C \alpha. \ee Now, taking (formally) ${\bf v} = \partial_t
(E_\alpha - E)$ in (\ref{eq2}) we arrive at
\[
\ds \eps_0\partial_t \int_{\Omega} \big[ |\partial_t (E_\alpha -
E)|^2 + \frac{1}{\mu_\alpha} |\curl (E_\alpha -E) |^2 \big]~dx =\]
\[ 2 \sum_{j=1}^m (\frac{1}{\mu_0} - \frac{1}{\mu_j})
\int_{z_j + \alpha B_j} \curl E \cdot \curl
\partial_t (E_\alpha -E)~dx.
\]

By using the regularity of $E$ in $\Omega$ and estimate
(\ref{estenergy1}) given above, we see that \bea \ds |
\sum_{j=1}^m (\frac{1}{\mu_0} - \frac{1}{\mu_j}) \int_{z_j +
\alpha B_j} \curl E \cdot \curl \partial_t (E_\alpha - E)~dx |&
\leq & C ||\curl E ||_{{\bf H}^2({\cal B}_\alpha)} ||  \partial_t
(E_\alpha - E)||_{{\bf H}^{-1}(\Omega)}\\
&\leq & C \alpha^2, \eea
 where $C$ is independent of $t$ and $\alpha$, and so, we
obtain
\[
\ds
\partial_t \int_{\Omega} \big[
|\partial_t (E_\alpha - E)|^2 + \frac{1}{\mu_\alpha} |\curl
(E_\alpha -E) |^2 \big]~dx \leq C \alpha^2
\]
which yields the following estimate
\[
|| \partial_t (E_\alpha - E) ||_{L^{\infty}(0, T; {\bf
L}^2(\Omega))} + || E_\alpha - E ||_{L^{\infty}(0, T;
X_N(\Omega))}  \leq C \alpha,
\]
where $C$ is independent of $\alpha$ and the points $\{
z_j\}_{j=1}^m$.

\square\\

Now, we can estimate $\curl E_\alpha-\curl E_0$ as follows.
\begin{proposition}\label{est-prop1}
Let $E_\alpha$ and $E$ be solutions to the problems
(\ref{alphwam}) and (\ref{wom}) respectively. There exist
constants $0<\alpha_0$, $C$ such that for $0<\alpha<\alpha_0$ the
following estimate holds: \be\label{estenergy2} \ds || \curl
(E_\alpha -E_0)||_{L^{\infty}(0, T; {\bf L}^2(\Omega))}\leq C
\alpha, \ee
\end{proposition}
\proof To prove estimate (\ref{estenergy2}) it is useful to
introduce the following function

\be\label{r2} \ds \hat{v}(x) = \int_0^T v(x, t) z(t) \; dt \in
L^2(\Omega), \ee where $v \in L^1(0, T; L^2(\Omega))$ and $z(t)$
is a given function in
 ${\cal C}^\infty_0(]0, T[)$.\\
Then,
\[
\ds \hat{E}(x) = \int_0^T E(x, t) z(t) \; dt {\rm \; and \;}
  \hat{E}_\alpha(x) = \int_0^T E_\alpha(x, t) z(t) \; dt
\in X(\Omega),
\]
which by relation (\ref{estenergy0}) give
\[
\left\{
\begin{array}{l}
(\hat{E}_\alpha
- \hat{E})  \in {\bf H}^1(\Omega),\\
\nm \ds \curl  \curl  (\hat{E}_\alpha - \hat{E}) =  0(\alpha)
\quad
{\rm in}\; \Omega,\\
\nm \dvg(\hat{E}_\alpha  - \hat{E})  = 0 \quad {\rm in}\;
\Omega,\\  \nm (\hat{E}_\alpha  - \hat{E}) \times {\bf n}
|_{\partial \Omega} = 0,
\end{array}
\right.
\]
and so, \be\label{r3} ||\curl(\hat{E}_\alpha - \hat{E})||_{{\bf
L}^2(\Omega)}=O(\alpha). \ee The fact that $\curl(E_\alpha - E)$
belongs to $L^\infty(0, T; {\bf L}^2(\Omega))$ and by using
estimate (\ref{r3}) we deduce that
\[
\int_{\Omega} |\curl E_\alpha(x,t) - \curl E(x,t)|^2~dx
=O(\alpha^2)\quad \mbox{ a.e. in }t\in (0,T),
\]
which means that
\[
||\curl (E_\alpha - E)||_{{\bf L}^2(\Omega)}=O(\alpha)\quad \mbox{
a.e. in }t\in (0,T).
\]
Thus, estimate (\ref{estenergy2}) follows immediately if we take
the sup on $t\in (0,T)$ in the last relation. \square\\

Before formulating our main result in this section, let us denote
$\Phi_j, j=1, \ldots, m$ the unique vector-valued solution to \be
\label{Phi} \left\{
\begin{array}{l} \Delta \Phi_j = 0 {\rm \; in\;} B_j, {\rm \; and \;}
\R^d \setminus \overline{B_j},\\ \nm \Phi_j {\rm \; is \;
continuous \; across\; } \partial B_j,\\ \nm \ds
\frac{\mu_j}{\mu_0} \frac{\partial \Phi_j}{\partial \nu_j}|_+ -
\frac{\partial \Phi_j}{\partial \nu_j}|_- = - \nu_j,\\ \nm \ds
\lim_{|y| \rightarrow + \infty} |\Phi_j(y)| = 0,
\end{array}
\right. \ee where $\nu_j$ denotes the outward unit normal to
$\partial B_j$, and superscripts $-$ and $+$ indicate the limiting
values as the point approaches $\partial B_j$ from outside $B_j$,
and from inside $B_j$, respectively. The existence and uniqueness
of this $\Phi_j$ can be established using single layer potentials
with suitably chosen densities, see \cite{CMV} for the case of
conductivity problem. For each inhomogeneity $z_j+ \alpha B_j$ we
introduce the  polarizability tensor $M_j$ which is a $d \times
d$, symmetric, positive definite matrix associated with the j-th
inhomogeneity, given by \begin{equation} \label{Meq} \ds (M_j)_{k,
l} = e_k \cdot (\int_{\partial B_j}
 (\nu_j +
(\frac{\mu_j}{\mu_0} - 1) \frac{\partial \Phi_j}{\partial
\nu_j}|_+ (y)) y \cdot e_l \; d\sigma_j(y)). \end{equation}
 Here $(e_1, \ldots, e_d)$ is an
orthonormal basis of $\R^d$. In terms of this function we are able
to prove the following result about the asymptotic behavior of
$\ds \curl E_\alpha \cdot\nu_j|_{\partial (z_j + \alpha B_j)^+}$.
\begin{theorem} \label{thm0} Suppose that (\ref{f1}) is satisfied and
 let $\Phi_j, j=1, \ldots, m$ be given as in (\ref{Phi}). Then, for the solutions
 $E_\alpha$, $E$ of problems (\ref{alphwam}) and (\ref{wom}) respectively, and for
  $y \in \partial
B_j$ we have \be \label{eq3} \ds(\curl E_\alpha (z_j + \alpha
y)\cdot\nu_j)|_{\partial (z_j + \alpha B_j)^+} = \curl E(z_j, t)
\cdot\nu_j \ee$$+ (1-\frac{\mu_j}{\mu_0})\frac{\partial
\Phi_j}{\partial \nu_j}|_+ (y) \cdot \curl E(z_j, t) +  o(1). $$
The term $o(1)$ uniform in $y \in \partial B_j$ and $t \in (0, T)$
and depends on the shape of $\{B_j \}_{j=1}^{m}$ and $\Omega$, the
constants $c_0$, $T$, $\mu_0$, $\{ \mu_j \}_{j=1}^{m}$, the data
$\varphi, \psi,$ and $f$, but is  otherwise independent of the
points $\{z_j \}_{j=1}^{m}$.
\end{theorem}
\proof

Let $\mathcal{E}_\alpha = \curl E_\alpha(x, t)$ and
$\mathcal{E}_0= \curl E(x, t)$. Then, according to
(\ref{alphwam})-(\ref{wom}) we have \be \label{u}
\eps_0\partial_t^2 E_\alpha - \curl \frac{1}{\mu_\alpha}
\mathcal{E}_\alpha = 0 {\rm \; and \;} \curl \mathcal{E}_\alpha =
0, {\rm \;  for \;}x\in\Omega. \ee We restrict, for simplicity,
our attention to the case of a single inhomogeneity, i.e., the
case $m=1$. The proof for any fixed number $m$ of well separated
inhomogeneities follows by iteration of the argument that we will
present for the case $m=1$. In order to further simplify notation,
we assume that the single inhomogeneity has the form $\alpha B$,
that is, we assume it is centered at the origin. We denote the
electromagnetic permeability inside $\alpha B$  by $\mu_*$ and
define $\Phi_*$ the same as $\Phi_j$, defined in (\ref{Phi}), but
with $B_j$ and $\mu_j$ replaced by $B$ and $\mu_*$, respectively.
Define $\nu$ to be the outward unit normal to $\partial B$. Now,
following a common practice in multiscale expansions
 we introduce the local variable $\ds y = \frac{x}{\alpha}$, then the domain
 $\tilde{\Omega} = \ds(\frac{\Omega}{\alpha})$ is well defined.\\
Next, let $\varpi$ be given in
 ${\cal C}^\infty_0(]0, T[)$. For any function
  $v \in {\bf L}^1(0, T; {\bf L}^2(\Omega))$, we define
\[
\hat{v}(x) = \ds  \int_0^T v(x, t) \,  \varpi(t) \; dt \in {\bf
L}^2(\Omega).
\]
We remark that $\widehat{\partial_t v} (x) = - \ds \int_0^T v(x,
t) \varpi^\prime(t) \; dt.$ So that we deduce from (\ref{u}) that
$\hat{\mathcal{E}}_\alpha$ satisfies
\[
\left\{ \begin{array}{l} \ds
 \curl \frac{1}{\mu_\alpha} \hat{\mathcal{E}}_\alpha = \int_0^T E_\alpha \,
\varpi^{\prime \prime}(t) \; dt
\quad {\rm in\;} \Omega,\\
\nm \curl \hat{\mathcal{E}}_\alpha = 0 \quad {\rm in\;} \Omega.
\end{array}
\right.
\]
Analogously,  $\hat{\mathcal{E}}$ satisfies
\[
\left\{ \begin{array}{l} \ds
 \frac{1}{\mu_0}\curl \hat{\mathcal{E}} = \int_0^T E  \, \varpi^{\prime \prime}(t) \; dt
\quad {\rm in\;} \Omega,\\
\nm \curl \hat{\mathcal{E}}  = 0 \quad {\rm in\;} \Omega.
\end{array}
\right.
\]
Indeed, we have $\hat{\mathcal{E}}_\alpha \times {\bf n} =
\hat{\mathcal{E}} \times {\bf n} = \curl_{\partial \Omega} \hat{f}
\times {\bf n}$ on the boundary $\partial \Omega$, where
$\curl_{\partial \Omega}$ is the tangential curl. Following
\cite{AVV} and \cite{A}, we introduce $q_\alpha^*$ as the unique
solution to the following problem
\[
\left\{ \begin{array}{l} \ds
 \Delta q_\alpha^* = 0
\quad {\rm in\;} \tilde{\Omega} = (\frac{\Omega}{\alpha})
\setminus \overline{B}
{\rm \; and \; in\;} B, \nonumber\\
\nm \ds
q_\alpha^* {\rm \;is\;continuous\;across\;} \partial B, \nonumber\\
\nm \ds \mu_0 \frac{\partial q_\alpha^*}{\partial \nu}|_+ - \mu_*
 \frac{\partial q_\alpha^*}{\partial \nu}|_{-} = - (\mu_0 - \mu_*)
\hat{\mathcal{E}}(\alpha y) \cdot \nu \quad {\rm on\;} \partial B,\\
\nm \ds q_\alpha^* = 0 \quad {\rm on\;} \partial \tilde{\Omega}.
\nonumber
\end{array}
\right.
\]
The jump condition
$$
\ds \mu_0 \frac{\partial q_\alpha^*}{\partial \nu}|_+ - \mu_*
 \frac{\partial q_\alpha^*}{\partial \nu}|_{-} = - (\mu_0 - \mu_*)
\hat{\mathcal{E}}(\alpha y) \cdot \nu \quad {\rm on\;} \partial B
$$
guarantees that $\hat{\mathcal{E}}_\alpha(x) -
\hat{\mathcal{E}}(x) - \grad_y q_\alpha^* (\frac{x}{\alpha})$
belongs to the functional space $X_N(\Omega)$, where
$\grad_{\partial \Omega}$ is the tangential gradient. Since
\[
\left\{ \begin{array}{l} \ds
 \curl \frac{1}{\mu_\alpha}( \hat{\mathcal{E}}_\alpha
- \hat{\mathcal{E}} - \grad_y q_\alpha^* (\frac{x}{\alpha}) ) =
\int_0^T \Bigr[ E_\alpha - \chi(\Omega \setminus \overline{\alpha
B}) E + \frac{\mu_*}{\mu_0} \chi(\alpha B) E \Bigr] \varpi^{\prime
\prime}(t) \; dt
\quad {\rm in\;} \Omega,\\
\nm \ds \curl ( \hat{\mathcal{E}}_\alpha - \hat{\mathcal{E}} -
\grad_y q_\alpha^* (\frac{x}{\alpha}))
= 0 \quad {\rm in\;} \Omega,\\
\nm ( \hat{\mathcal{E}}_\alpha - \hat{\mathcal{E}} - \grad_y
q_\alpha^* (\frac{x}{\alpha}) ) \times {\bf n} = 0 \quad {\rm
on\;}
\partial \Omega,
\end{array}
\right.
\]
where $\chi(\omega)$ is the characteristic function of the domain
$\omega$, we arrive, as a consequence of the energy estimate given
by Lemma \ref{est-lem1}, at the following
\[
\left\{ \begin{array}{l}
 (\hat{\mathcal{E}}_\alpha
- \hat{\mathcal{E}} - \grad_y q_\alpha^* (\frac{x}{\alpha}))
\in X_N(\Omega),\\
\nm \ds
 \curl \frac{1}{\mu_\alpha} ( \hat{\mathcal{E}}_\alpha
- \hat{\mathcal{E}} - \curl_y q_\alpha^* (\frac{x}{\alpha})) =
0(\alpha)
\quad {\rm in\;} \Omega,\\
\nm \ds \curl ( \hat{\mathcal{E}}_\alpha - \hat{\mathcal{E}} -
\grad_y q_\alpha^* (\frac{x}{\alpha}))
= 0 \quad {\rm in\;} \Omega,\\
\nm ( \hat{\mathcal{E}}_\alpha - \hat{\mathcal{E}} - \grad_y
q_\alpha^* (\frac{x}{\alpha}))  \times {\bf n} = 0 \quad {\rm
on\;}
\partial \Omega.
\end{array}
\right.
\]
From \cite{AVV} we know that this yields the following estimate
\[
|| \curl \frac{1}{\mu_\alpha}( \hat{\mathcal{E}}_\alpha -
\hat{\mathcal{E}} - \grad_y q_\alpha^* (\frac{x}{\alpha}) )
||_{L^2(\Omega)} + ||  \hat{\mathcal{E}}_\alpha -
\hat{\mathcal{E}} - \grad_y q_\alpha^* (\frac{x}{\alpha})
||_{L^2(\Omega)} \leq C \alpha,
\]
and so,
\[
( \hat{\mathcal{E}}_\alpha - \hat{\mathcal{E}} - \grad_y
q_\alpha^* (\frac{x}{\alpha})) \cdot \nu |_+ = 0(\alpha) \quad
{\rm on\;}
\partial (\alpha B).
\]
Now, we denote by $q_*$ be the unique (scalar) solution to
\[
\left\{ \begin{array}{l} \ds
 \Delta q_* = 0
\quad {\rm in\;} \RR^d \setminus \overline{B}
{\rm \; and \; in\;} B,\\
\nm \ds
q_* {\rm \;is\;continuous\;across\;} \partial B, \\
\nm \ds \mu_0 \frac{\partial q_*}{\partial \nu}|_+ - \mu_*
 \frac{\partial q_*}{\partial \nu}|_{-} = - (\mu_0 - \mu_*)
\hat{\mathcal{E}}(0) \cdot \nu \quad {\rm on\;} \partial B,\\
\nm \ds \lim_{|y| \rightarrow + \infty} q_* = 0.
\end{array}
\right.
\]
In the spirit of Theorem 1 in \cite{CMV} it follows that
\[
\ds ||(\grad_y q_* - \grad_y q_\alpha^*)(\frac{x}{\alpha})
||_{{\bf L}^2(\Omega)} \leq C \alpha^{1/2},
\]
which yields
\[
( \hat{\mathcal{E}}_\alpha - \hat{\mathcal{E}} - \grad_y q_*
(\frac{x}{\alpha})) \cdot \nu = o(1) \quad {\rm on\;} \partial
(\alpha B).
\]
Writing $q_*$ in terms of $\Phi_*$ gives
\[
\ds \int_0^T \Bigr[ (\curl E_\alpha(\alpha y)\cdot\nu)|_{\partial
(\alpha B)^+}  -  \nu \cdot \curl E(0, t) - (\frac{\mu_0}{\mu_*} -
1) \frac{\partial \Phi_*}{\partial \nu}|_+ (y) \cdot \curl E(0, t)
\Bigr] \varpi(t) \; dt =  o(1),
\]
for any $\varpi \in {\cal C}^\infty_0(]0, T[)$, and so, by
iterating the same argument for the case of $m$ (well separated)
inhomogeneities $z_j + \alpha B_j, j=1, \ldots, m$, we arrive at
the promised asymptotic formula (\ref{eq3}).

\square

\section{The identification procedure} Before describing our
identification procedure, let us introduce the following cutoff
function $\beta(x)  \in {\cal C}^{\infty}_0(\Omega)$ such that
$\beta \equiv 1$ in a subdomain $\Omega^{\prime}$ of $\Omega$ that
contains the inhomogeneities ${\cal B}_\alpha$ and let $\eta \in
\RR^d$. We will take in what follows $E(x, t) = \eta^{\perp} e^{i
\eta \cdot x -i\sqrt{\mu_0}|\eta | t}$ where $\eta^{\perp}$ is a
unit vector that is orthogonal to $\eta$ which corresponds to
taking $\varphi(x) = \eta^{\perp} e^{i \eta \cdot x}, \psi(x) = -
i \sqrt{\mu_0}|\eta| \eta^{\perp} e^{i \eta \cdot x}, $ and $f(x,
t) = \eta^{\perp}\times {\bf n} e^{i \eta\cdot x - i
\sqrt{\mu_0}|\eta| t}$ and assume that we are in possession of the
measurements of:
$$\curl E_\alpha \times {\bf n} \quad \mbox{ on } \Gamma \times (0, T),$$
 where $\Gamma$ is an open part of $\partial \Omega$.  Suppose now that $T$ and
the part $\Gamma$ of the boundary $\partial \Omega$ are such that
they  geometrically control $\Omega$ which roughly means that
every geometrical optic ray, starting at any point $x \in \Omega$
at time $t=0$ hits $\Gamma$ before time $T$ at a non diffractive
point, see \cite{BLR}. It follows from \cite{NI} (see also
\cite{LAG}, \cite{K} and \cite{KO}) that there exists (a unique)
$g_\eta \in H^1_0(0, T; TL^2(\Gamma))$ (constructed by the Hilbert
Uniqueness Method) such that the unique weak solution $w_\eta$ to
\be \label{wetam} \left\{
\begin{array}{l}
\ds (\partial_t^2 + \curl  \curl ) w_\eta  = 0 \quad
{\rm in}\; \Omega \times (0, T),\\
\nm \dvg w_\eta = 0 \quad
{\rm in}\; \Omega \times (0, T),\\
\nm w_\eta  |_{t=0} = \beta(x) \eta^{\perp} e^{i \eta \cdot x},
\partial_t w_\eta  |_{t=0} = 0 \quad {\rm in}\; \Omega,\\ \nm
w_\eta  \times {\bf n} |_{\partial \Omega \setminus
\overline{\Gamma}
 \times (0, T)} = 0,\\
\nm w_\eta  \times {\bf n} |_{\Gamma
 \times (0, T)} = g_\eta,
\end{array}
\right. \ee satisfies $w_\eta(T) = \partial_t w_\eta (T) = 0$ in
$\Omega$.

Let $\theta_\eta \in H^1(0, T; TL^2(\Gamma))$ denote the unique
solution of the Volterra equation of second kind \be\label{eq4m}
\left\{\begin{array}{l} \ds
\partial_t \theta_\eta (x, t) + \int_t^T e^{- i | \eta | (s -t)}
( \theta_\eta (x, s) - i | \eta | \partial_t \theta_\eta (x, s))
\; ds = g_\eta(x, t) \quad {\rm for\;} x \in \Gamma, t \in (0, T),
\\ \theta_\eta(x, 0) = 0 \quad {\rm for\;} x \in \Gamma.
\end{array}
\right. \ee The existence and uniqueness of this $\theta_\eta$ in
$H^1(0, T; TL^2(\Gamma))$ for any $\eta \in \RR^d$ can be
established using the resolvent kernel. However, observing from
differentiation of (\ref{eq4m}) with respect to $t$ that
$\theta_\eta$ is the unique solution of the ODE: \be \label{eq4p}
\left\{
\begin{array}{l}
\ds
\partial_t^2 \theta_\eta - \theta_\eta = e^{i |\eta| t}
\partial_t ( e^{-i |\eta| t} g_\eta) \quad {\rm for\;} x \in
\Gamma, t \in (0, T), \\ \theta_\eta(x, 0) = 0,
\partial_t \theta_\eta(x, T) = 0  \quad {\rm for\;} x \in \Gamma,
\end{array}
\right. \ee the function $\theta_\eta$ may be find (in practice)
explicitly with variation of parameters and it also immediately
follows from this observation that $\theta_\eta$ belongs to
$H^2(0, T; TL^2(\Gamma))$.\\
We introduce $v_{ \eta}$ as the unique weak solution (obtained by
transposition as done in \cite{LM} and in \cite{L} [Theorem 4.2,
page 46] for the scalar function) in $ {\cal C}^0(0, T; X(\Omega))
\cap {\cal C}^1(0, T; L^2(\Omega))$ to the following problem
\[
\left\{
\begin{array}{l}
\ds (\partial_t^2 + \curl  \curl ) v_\eta  = 0 \quad
{\rm in}\; \Omega \times (0, T),\\
\nm \dvg v_\eta = 0 \quad
{\rm in}\; \Omega \times (0, T),\\
\nm v_\eta  |_{t=0} = 0  \quad {\rm in}\; \Omega,\\
\nm \ds
\partial_t v_\eta  |_{t=0} = \sum_{j=1}^m i (1 - \frac{\mu_0}{\mu_j}) \eta \times
 (\nu_j +
(\frac{\mu_0}{\mu_j} - 1) \frac{\partial \Phi_j}{\partial
\nu_j}|_+) e^{i \eta \cdot z_j} \delta_{\partial (z_j + \alpha
B_j)}  \in Y(\Omega) \quad {\rm in}\; \Omega,\\ \nm v_\eta  \times
{\bf n} |_{\partial \Omega
 \times (0, T)} = 0.
\end{array}
\right.
\]

Then, the following holds.

\begin{proposition}\label{p4.1}
Suppose that $\Gamma$ and $T$ geometrically control $\Omega$. For
any $\eta \in \RR^d$ we have \be\label{rel-prop3.1} \int_0^T
\int_\Gamma g_\eta \cdot(\curl v_{\eta}\times{\bf
n})~d\sigma(x)dt=\alpha^2 \sum_{j=1}^m  \mu_0(1 -
\frac{\mu_j}{\mu_0}) e^{2 i \eta \cdot z_j} \eta \cdot
\int_{\partial B_j}
 (\nu_j \ee \[+
(\frac{\mu_j}{\mu_0} - 1) \frac{\partial \Phi_j}{\partial
\nu_j}|_+ (y)) \eta \cdot y \; ds_j(y)+ o(\alpha^2).
\]
\end{proposition}
\proof Multiply the equation $\ds (\partial_t^2 + \curl  \curl )
v_\eta  = 0$ by $w_\eta$ and integrating by parts over $(0, T)
\times \Omega$, for any $\eta \in \RR^d$ we have
\[
\alpha \sum_{j=1}^m i (1 - \frac{\mu_j}{\mu_0}) e^{2 i \eta \cdot
z_j} \eta \cdot \int_{\partial B_j}
 (\nu_j +
(\frac{\mu_j}{\mu_0} - 1) \frac{\partial \Phi_j}{\partial
\nu_j}|_+ (y)) e^{i \alpha \eta \cdot y} \; ds(y)  = \]
\[-\mu_{0}^{-1} \int_0^T \int_\Gamma g_\eta \cdot(\curl v_{
\eta}\times{\bf n})~d\sigma(x)dt.
\]
Now, we take the Taylor expansion of $\alpha e^{i \alpha \eta
\cdot y}$ in the left side of the last equation, we obtain the
convenient asymptotic formula (\ref{rel-prop3.1}).

\square

 To identify the locations and certain properties of
the small inhomogeneities ${\cal B}_\alpha$ let us view the
averaging of the boundary measurements \[\ds \curl E_\alpha
\times{ \bf n} |_{\Gamma \times (0, T)},\] using the solution
$\theta_\eta$ to the Volterra equation (\ref{eq4m}) or
equivalently the ODE (\ref{eq4p}), as a function of $\eta$. The
following holds.\\

\begin{theorem} \label{th2}
Let $\eta \in \RR^d$. Let $E_\alpha$ be the unique solution in $
{\cal C}^0(0, T; X(\Omega)) \cap {\cal C}^1(0, T; L^2(\Omega))$ to
the Maxwell's equations (\ref{alphwam}) with $ \varphi(x) =
\eta^{\perp} e^{i \eta \cdot x},$ $ \psi(x) = - i \sqrt{\mu_0}|
\eta | \eta^{\perp} e^{i \eta \cdot x},$ and $ f(x, t) =
\eta^{\perp} e^{i \eta \cdot x  -  i \sqrt{\mu_0}| \eta | t}.$
Suppose that $\Gamma$ and $T$ geometrically control $\Omega$, then
we have \be \label{eqam}
\begin{array}{l}
\ds \int_0^T \int_\Gamma \Big[ \theta_\eta \cdot (\curl E_\alpha
\times {\bf n} - \curl E \times {\bf n}) +
\partial_t \theta_\eta \cdot \partial_t (\curl E_\alpha \times {\bf n}
- \curl E \times {\bf n})  \Big]~d\sigma(x)dt =  \\ \nm \ds
\alpha^2 \sum_{j=1}^m (\mu_0 - \mu_j) e^{2 i \eta \cdot z_j}
M_j(\eta) \cdot \eta \;
  + O(\alpha^2),
  \end{array}
\ee where $\theta_\eta$ is the unique solution to the Volterra
equation (\ref{eq4p}) with $g_\eta$ defined as the boundary
control  in (\ref{wetam}) and $M_j$ is the polarization tensor of
$B_j$, defined by \be \label{pt} \ds (M_j)_{k, l} = e_k \cdot
(\int_{\partial B_j}
 (\nu_j +
(\frac{\mu_j}{\mu_0} - 1) \frac{\partial \Phi_j}{\partial
\nu_j}|_+ (y)) y \cdot e_l \; ds_j(y)). \ee Here $(e_1, e_2)$ is
an orthonormal basis of $\RR^d$. The term $O(\alpha^2)$ is
independent of the points $\{z_j,\quad j=1,\cdots, m\}$.
\end{theorem}
\proof From $\partial_t \theta_\eta (T) = 0$ and $(\curl E_\alpha
\times {\bf n} - \curl E \times {\bf n})|_{t=0}  = 0$ the term $
\ds \int_0^T \int_\Gamma
\partial_t \theta_\eta \cdot\partial_t (\curl E_\alpha
\times {\bf n} - \curl E \times {\bf n})~d\sigma(x)dt$ has to be
interpreted as follows \be \label{d} \ds \int_0^T \int_\Gamma
\partial_t \theta_\eta\cdot \partial_t (\curl E_\alpha
\times {\bf n} - \curl E \times {\bf n})  = - \int_0^T \int_\Gamma
\partial^2_t \theta_\eta \cdot(\curl E_\alpha
\times {\bf n} - \curl E \times {\bf n}). \ee Next, introduce
\be\label{def1} \ds \tilde{E}_{\alpha, \eta}(x, t) = E(x, t) +
\int_0^t e^{- i \sqrt{\mu_0} | \eta| s} v_{\eta}(x, t-s)\; ds, x
\in \Omega, t \in (0, T). \ee We have
\[\begin{array}{l}
\ds \int_0^T \int_\Gamma \Bigr[ \theta_\eta \cdot(\curl E_\alpha
\times {\bf n} - \curl E \times {\bf n} ) +
\partial_t \theta_\eta\cdot \partial_t (\curl E_\alpha
\times {\bf n} - \curl E \times {\bf n} ) \Bigr]
= \\
\nm \ds \int_0^T \int_\Gamma \Bigr[ \theta_\eta \cdot(\curl
E_\alpha \times {\bf n}- \curl \tilde{E}_{\alpha,\eta}\times {\bf
n}) +
\partial_t \theta_\eta\cdot \partial_t (\curl
E_\alpha \times {\bf n}- \curl \tilde{E}_{\alpha,\eta}\times {\bf
n} ) \Bigr]\\
\nm \ds + \int_0^T \int_{\Gamma} \Bigr[ \theta_\eta\cdot \int_0^t
e^{- i \sqrt{\mu_0} | \eta| s} v_{\eta}(x, t-s)\times{\bf n}\; ds
+
\partial_t \theta_\eta \cdot\partial_t \int_0^t e^{- i
\sqrt{\mu_0} | \eta| s} v_{\eta}(x, t-s)\times{\bf n}\; ds \Bigr].
\end{array}
\]
Since $\theta_\eta$ satisfies the Volterra equation (\ref{eq4p})
and
\[\begin{array}{l}
\ds \partial_t ( \int_0^t e^{-  i \sqrt{\mu_0} | \eta | s}
v_{\eta}(x, t-s)\times{\bf n}\; ds ) =
\partial_t (- e^{- i \sqrt{\mu_0} |\eta| t}  \int_0^t e^{ i \sqrt{\mu_0}
| \eta | s} v_{\eta}(x, s)\times{\bf n}\; ds )
\\ \nm \ds =  i \sqrt{\mu_0} |\eta| e^{- i \sqrt{\mu_0} |\eta| t}
 \int_0^t e^{ i \sqrt{\mu_0}| \eta | s}
v_{\eta}(x, s)\times{\bf n}\; ds + v_{\eta}(x, t)\times{\bf n},
\end{array}
\]
we obtain by integrating by parts over $(0, T)$ that
\[\begin{array}{l}
  \ds \int_0^T \int_{\Gamma} \Bigr[
\theta_\eta\cdot \int_0^t e^{-  i\sqrt{\mu_0}  | \eta | s}
v_{\eta}(x, t-s)\times{\bf n}\; ds +
\partial_t \theta_\eta\cdot
\partial_t \int_0^t e^{- i\sqrt{\mu_0} | \eta | s} v_{\eta}(x, t-s)\times{\bf n}\;
 ds \Bigr]  \\= \ds
\int_0^T \int_{\Gamma} ( v_{\eta}(x, t)\times{\bf n})\cdot
(\partial_t \theta_\eta + \int_t^T \theta_\eta(s) e^{i
\sqrt{\mu_0} |\eta| (t-s)} \;ds)\\ \nm \quad \ds  - i \sqrt{\mu_0}
|\eta| (e^{- i \sqrt{\mu_0}|\eta| t}
\partial_t \theta_\eta (t))\cdot \int_0^t e^{ i\sqrt{\mu_0} | \eta | s}
v_{\eta}(x, s)\times{\bf n}\; ds\; dt\\ \nm = \ds \int_0^T
\int_{\Gamma} v_{\eta}(x, t)\times{\bf n}\cdot (\partial_t
\theta_\eta + \int_t^T (\theta_\eta(s) - i \sqrt{\mu_0} |\eta|
\partial_t \theta_\eta(s))
 e^{ i \sqrt{\mu_0}| \eta |(t- s)} \;ds) \;dt\\
\nm \ds
 = \int_0^T \int_{\Gamma} g_\eta(x, t)
\cdot (\curl v_{\eta}(x, t)\times{\bf n})\; dt
\end{array}
\]
and so, from Proposition \ref{p4.1} we obtain
\[\begin{array}{l}
\ds \int_0^T \int_\Gamma
 \Bigr[
\theta_\eta \cdot(\curl E_\alpha \times {\bf n} - \curl E \times
{\bf n}  ) +
\partial_t \theta_\eta \cdot\partial_t (\curl E_\alpha
\times {\bf n} - \curl E \times {\bf n} ) \Bigr]
= \\
\nm \ds \alpha^2  \sum_{j=1}^m  (1 - \frac{\mu_j}{\mu_0}) e^{2 i
\eta \cdot z_j} \eta \cdot \int_{\partial B_j}
 (\nu_j +
(\frac{\mu_j}{\mu_0} - 1) \frac{\partial \Phi_j}{\partial
\nu_j}|_+ (y)) \eta \cdot y \; ds_j(y)\\
\nm \ds + \int_0^T \int_\Gamma
 \Bigr[
\theta_\eta \cdot(\curl E_\alpha \times {\bf n} - \curl
\tilde{E}_{\alpha, \eta} \times {\bf n} ) +
\partial_t \theta_\eta \cdot\partial_t (\curl E_\alpha \times {\bf n}\\
\nm \ds  - \curl \tilde{E}_{\alpha, \eta} \times {\bf n} ) \Bigr]
+ o(\alpha^2).
\end{array}
\]
In order to prove Theorem {\ref{th2} it suffices then to show that
\be \label{p}
 \int_0^T \int_\Gamma
 \Bigr[
\theta_\eta \cdot(\curl E_\alpha \times {\bf n} - \curl
\tilde{E}_{\alpha, \eta} \times {\bf n} ) +
\partial_t \theta_\eta\cdot \partial_t (\curl E_\alpha \times {\bf n} - \curl
\tilde{E}_{\alpha, \eta} \times {\bf n} ) \Bigr] = o(\alpha^2).
\ee Since
\[
\left\{
\begin{array}{l}
\ds (\partial_t^2 - \curl\frac{1}{\mu_0}\curl) ( \int_0^t e^{- i
\sqrt{\mu_0} | \eta| s} v_{ \eta}(x, t-s)\; ds)
\\ \ds = \sum_{j=1}^m i (1 - \frac{\mu_j}{\mu_0}) \eta \times
(\nu_j + (\frac{\mu_j}{\mu_0} - 1) \frac{\partial \Phi_j}{\partial
\nu_j}|_+ (y)) e^{i \eta \cdot z_j} \delta_{\partial (z_j + \alpha
B_j)} e^{- i \sqrt{\mu_0} |\eta| t} \quad {\rm in}\; \Omega \times
(0, T),\\ \nm \ds ( \int_0^t e^{- i \sqrt{\mu_0} | \eta| s} v_{
 \eta}(x, t-s)\; ds) |_{t=0} = 0,
\partial_t ( \int_0^t e^{- i \sqrt{\mu_0} | \eta| s} v_{ \eta}(x, t-s)\; ds)|_{t=0} = 0
 \quad {\rm in}\; \Omega,\\ \nm
\ds ( \int_0^t e^{- i \sqrt{\mu_0} | \eta| s} v_{\eta}(x, t-s)\;
ds)\times{\bf n} |_{\partial \Omega \times (0, T)} = 0,
\end{array}
\right.
\]
it follows from Theorem \ref{thm0} that
\[
\left\{
\begin{array}{l}
\ds (\partial_t^2 - \curl\frac{1}{\mu_0}\curl) (E_\alpha -
\tilde{E}_{\alpha, \eta})  = o(\alpha^2) \quad {\rm in}\; \Omega
\times (0, T),\\ \nm (E_\alpha - \tilde{E}_{\alpha, \eta}) |_{t=0}
= 0,
\partial_t (E_\alpha - \tilde{E}_{\alpha, \eta}) |_{t=0} = 0
 \quad {\rm in}\; \Omega,\\ \nm
(E_\alpha - \tilde{E}_{\alpha, \eta})\times{\bf n}|_{\partial
\Omega \times (0, T)} = 0.
\end{array}
\right.
\]
Following the proof of Proposition \ref{est-prop1}, we immediately
obtain
\[
|| E_\alpha - \tilde{E}_{\alpha, \eta} ||_{L^2(\Omega)} =
o(\alpha^2), t \in (0, T), x \in \Omega,
\]
where $o(\alpha^2)$ is independent of the points $\{
z_j\}_{j=1}^m$. To prove (\ref{p}) it suffices then from (\ref{d})
to show that the  following estimate holds
\[
\ds || \curl E_\alpha \times {\bf n} - \curl \tilde{E}_{\alpha,
\eta} \times {\bf n}||_{L^2(0, T; TL^2(\Gamma))} = o(\alpha^2).
\]
Let
 $\theta$ be given in
 ${\cal C}^\infty_0(]0, T[)$ and define
\[
\ds\hat{\tilde{E}}_{\alpha, \eta}(x) = \ds  \int_0^T
\tilde{E}_{\alpha, \eta}(x, t) \theta(t) \; dt
\]
and
\[
\ds\hat{E}_{\alpha}(x) = \ds  \int_0^T E_{\alpha}(x, t) \theta(t)
\; dt.
\]
From definition (\ref{def1}) we can write \be \label{r8} \left\{
\begin{array}{l}
(\hat{E}_\alpha
- \hat{\tilde{E}}_\alpha)  \in {\bf H}^1(\Omega),\\
\nm \ds \curl  \curl  (\hat{E}_\alpha - \hat{\tilde{E}}_\alpha) =
0(\alpha) \in Y(\Omega) \quad
{\rm in}\; \Omega,\\
\nm \dvg(\hat{E}_\alpha  - \hat{\tilde{E}}_\alpha)  = 0 \quad {\rm
in}\; \Omega,\\  \nm (\hat{E}_\alpha  - \hat{\tilde{E}}_\alpha)
\times {\bf n} |_{\partial \Omega} = 0.
\end{array}
\right. \ee In the spirit of the standard elliptic regularity
\cite{E} we deduce for the boundary value problem (\ref{r8}) that
\[\ds || \curl(\hat{E}_\alpha - \hat{\tilde{E}}_{\alpha})\times {\bf n}
||_{{\bf L}^2(\Gamma)} = O(\alpha^2), \] for all $\theta \in {\cal
C}^\infty_0(]0, T[)$; whence
$$
||\curl(E_\alpha - \hat{E}_{\alpha})\times {\bf n}||_{{\bf
L}^2(\Gamma)}
 = o (\alpha^2) {\rm \; a.\, e.\; in \;} t \in (0, T), $$
and so, the desired estimate (\ref{eqam}) holds. The proof of
Theorem \ref{th2} is
then over. \square\\

Our identification procedure is deeply based on Theorem \ref{th2}.
 Let us neglect the asymptotically small remainder in the asymptotic formula
(\ref{eqam}), and define $\aleph_\alpha(\eta)$ by
\[
\ds \aleph_\alpha(\eta) =  \int_0^T \int_\Gamma
 \Bigr[
\theta_\eta \cdot(\curl(E_\alpha - E)\times {\bf n}) +
\partial_t \theta_\eta \cdot\partial_t (\curl(E_\alpha - E)\times {\bf n} ) \Bigr].
\]
Recall that the function $e^{2 i \eta \cdot z_j}$ is exactly the
Fourier Transform (up to a multiplicative constant) of the Dirac
function $\delta_{-2 z_j}$ (a point mass located at $- 2 z_j$).
From Theorem \ref{th2} it follows that the function $e^{2 i \eta
\cdot z_j}$ is (approximately) the Fourier Transform of a linear
combination of derivatives of point masses, or
\[
 \breve{\aleph}_\alpha(\eta) \approx \alpha^2 \sum_{j=1}^m L_j \delta_{-2 z_j},
\]
where $L_j$ is a second order constant coefficient, differential
operator whose
 coefficients depend on the polarization tensor $M_j$ defined by (\ref{pt})
(see  \cite{CMV} for its properties) and $
\breve{\aleph}_\alpha(\eta)$ represents the inverse Fourier
Transform of $\aleph_\alpha(\eta)$. The reader is referred to
\cite{CMV} for properties of the tensor polarization $M_j$.

The method of reconstruction consists in sampling values of
  $ \breve{\aleph}_\alpha(\eta)$ at some discrete set of points and then calculating
  the corresponding discrete inverse Fourier Transform. After a rescaling
  the support of this discrete inverse Fourier Transform yields the
location of the small inhomogeneities ${\cal B}_\alpha$. Once the
locations are known we may calculate the polarization tensors
$(M_j)^m_{j=1}$ by solving an appropriate linear system arising
from (\ref{eqam}).  This procedure generalizes the approach
developed in \cite{AMV} for the two-dimensional (time-independent)
inverse conductivity problem and generalize the results in
\cite{A} to the full time-dependent Maxwell's equations.
\section{Conclusion }
In this paper, we are convinced that the use of approximate
formulae such as (\ref{eqam}) represents a very promising approach
to the dynamical identification of small inhomogeneities that are
embedded in a homogeneous medium. We also believe that our method
yields a good approximation to small amplitude perturbations in
the electromagnetic parameters (for the example of electric
permittivity $\varepsilon_\alpha(x) = \varepsilon_0 + \alpha
\varepsilon(x)$) from the measurements:
$$\ds \curl
H_\alpha\times{\bf n}\quad \mbox{ on } \Gamma \times (0, T).$$ Our
method may yield the Fourier transform of the amplitude
perturbation $\varepsilon(x)$. This issue will be considered in a
forthcoming work \cite{khelifi-daveau2}.


\end{document}